\documentclass[12pt]{article}

\newtheorem{theorem}{Theorem}
\newtheorem{lemma}{Lemma}
\newtheorem{proposition}{Proposition}

\usepackage{amssymb}
\usepackage{amsmath}
\usepackage{latexsym}

\newcommand{\C}{{\mathbb C}}

\newcommand{\Z}{{\mathbb Z}}

\newcommand{\A}{\mathcal{A}}

\newcommand{\LL}{\mathcal{L}}

\begin{document}

\title{The $q$-Onsager algebra}

\author{
Tatsuro Ito 
and
Paul Terwilliger
}

\date{}
\maketitle

This article gives a summary of the finite-dimesional 
irreducible representations of the $q$-Onsager algebra, 
which are treated in detail in \cite{ITf}. 
 
\medskip
Fix a nonzero scalar $q \in \C$ 
which is not a root of unity. 
Let $\A=\A_q$ 
denote the associative $\C$-algebra with $1$ defined by 
genarators $z, z^*$ subject to the relations 
\[({\rm TD})\left\{
\begin{array}{lll}
[z,~z^2z^*-\beta z\,z^*z+z^*z^2] 
&=& \delta \,[z,~z^*],
\\
\lbrack z^*,~{z^*}^2 z -\beta z^*z\,z^* +z \,{z^*}^2 \rbrack
&=&\delta \,[z^*,~z], 
\end{array}
\right.
\]
where $\beta = q^2 + q^{-2}$ and $\delta = -{(q^2 - q^{-2})}^2$. 
$({\rm TD})$ can be regarded as a $q$-analogue of 
the Dolan- Grady relations 
and we call $\A$ the \textit{q-Onsager algebra}. 
We classify the finite-dimensional 
irreducible representations of $\A$. 
All such representations are explicitly constructed via embeddings 
of $\A$ into the $U_q(sl_2)$-loop algebra.
As an application, tridiagonal pairs of $q$-Racah type 
over $\C$ are classified 
in the case where $q$ is not a root of unity.

\medskip
The $U_q(sl_2)$-loop algebra 
${\LL}=U_q(L(sl_2))$ is the associative 
$\C$-algebra with 1 generated by 
$e^+_i, e^-_i, k_i, k^{-1}_i ~(i=0,1)$ subject to the relations 
\begin{eqnarray*}
k_0 k_1 &=& k_1 k_0 = 1, \\
k_i k^{-1}_i &=& k^{-1}_i k_i = 1, \\
k_i e^{\pm}_i k^{-1}_i &=& q^{\pm 2} e^{\pm}_i, \\
k_i  e^{\pm}_j k_i^{-1} &=& q^{\mp 2} e^{\pm}_j ~~~~(i \neq j), \\
\lbrack e^+_i, e^-_i \rbrack &=& \frac{{k_i} - k^{-1}_i}{q-q^{-1}}, \\
\lbrack e^+_i,e^-_j\rbrack&=&0 ~~~~(i \neq j), 
\end{eqnarray*}
\vspace{-9mm}
\begin{eqnarray*}
\lbrack e^{\pm}_i, (e^{\pm}_i)^2 e^{\pm}_j - (q^2 + q^{-2})
e^{\pm}_i e^{\pm}_j e^{\pm}_i + e^{\pm}_j (e^{\pm}_i )^2\rbrack=0 ~~~~
(i \neq j).
\end{eqnarray*}
Note that if we replace $k_0 k_1 = k_1 k_0 = 1$ 
in the defining relations for $\LL$ by $k_0 k_1 = k_1 k_0$, 
then we have the quantum affine algebra $U_q(\widehat{sl}_2)$: $\LL$ 
is isomorphic to the quotient algebra of $U_q(\widehat{sl}_2)$ 
by the two-sided ideal generated by $k_0 k_1 -1$.

\begin{proposition}
\label{prop: phi}
For arbitrary nonzero $s, t \in \C$, there exists an 
algebra homomorphism $\varphi_{s,t}$ from $\A$ to $\LL$ that sends 
$z,\, z^*$ to 
\begin{eqnarray*}
z_t(s) & =& x(s)+t\,k(s)+t^{-1} k(s)^{-1},\\
z^*_t(s) & =& y(s)+t^{-1} k(s) + t\, k(s)^{-1}, 
\end{eqnarray*}
respectively, where 
\begin{eqnarray*}
x(s) & =& \alpha(se_0^+ + s^{-1} e_1^- k_1) 
~~~ {\rm with} ~\alpha = - q^{-1} {(q-q^{-1})}^2,\\
y(s) & =& s e_0^- k_0 + s^{-1} e_1^+,\\
k(s) & =& s k_0. 
\end{eqnarray*}
Moreover $\varphi_{s,t}$ is injective. 
\end{proposition}

\medskip
We give an overview of finite-dimensional representations of 
$\LL$ that we need to state our explicit construction of 
irreducible $\A$-modules via $\varphi_{s,t}$. 
For $a \in \C ~(a \neq 0) $ and $\ell \in \Z ~(\ell > 0)$, 
$V(\ell,a)$ denotes the {\it evaluation module} 
of $\LL$, i.e., $V(\ell,a)$ is an $(\ell + 1)$-dimensional 
vector space over $\C$ with a basis 
$v_0, v_1, \ldots , v_{\ell}$ on which $\LL$ acts as follows: 
\begin{eqnarray*}
k_0 v_i  &=& q^{2i-\ell} \,v_i ,\\
k_1 v_i  &=& q^{\ell-2i} v_i, \\
e^+_0 v_i  &=& a\, q \,[i+1] \,v_{i+1}, \\
e^-_0 v_i  &=& a^{-1} q^{-1} [\ell-i+1]\, v_{i-1}, \\
e^+_1 v_i  &=&  [\ell-i+1] \,v_{i-1}, \\
e^-_1 v_i  &=& [i+1]\, v_{i+1},
\end{eqnarray*}
where $v_{-1} = v_{\ell+1} = 0 $ and 
$ [j] = [j]_q=(q^j - q^{-j})/(q - q^{-1})$. 
$V(\ell,a)$  is an irreducible  $\LL$-module. 
We call $v_0, v_1, \ldots , v_{\ell}$ a {\it standard basis}.

\medskip
Let $\Delta$ denote the {\it coproduct} of
$\LL$: the algebra homomorphism from $\LL$ to 
$\LL \otimes \LL$ defined by  
\begin{eqnarray*}
\Delta (k^{\pm 1}_i) & =& k^{\pm 1}_i \otimes k^{\pm 1}_i ,\\
\Delta (e^+_i) & =& k_i \otimes e^+_i + e^+_i \otimes 1, \\
\Delta (e^-_i k_i) & =& k_i \otimes e^-_i k_i + e^-_i k_i \otimes 1.
\end{eqnarray*}
Given $\LL$-modules $V_1,V_2$, the tensor product $V_1 \otimes V_2$ 
becomes an $\LL$-module via $\Delta$.
Given a set of evaluation modules $V(\ell_i,a_i)~(1 \leq i \leq n)$ 
of $\LL$, the tensor product 
$$ V(\ell_1,a_1) \otimes \cdots \otimes V(\ell_n,a_n) $$
makes sense as an $\LL$-module without being affected by   
the parentheses for the tensor product 
because of the coassociativity of $\Delta$. 

\medskip
With an evaluation module $V(\ell, a)$ of $\LL$, 
we associate the set $S(\ell, a)$ of scalars 
$a\,q^{-\ell+1}$, $a\,q^{-\ell+3}$, $\cdots$, $a\,q^{\ell-1}$: 
$$ S(\ell,a) = \{ a\,q^{2i-\ell+1} ~|~ 0 \leq i \leq \ell-1 \}. $$
The set $S(\ell, a)$ is called a $q$-$string$ of length $\ell$.
Two $q$-strings $S(\ell, a)$, $S(\ell', a')$ are said to be 
{\it adjacent} 
if $S(\ell, a) \cup S(\ell', a')$ is a longer $q$-string, i.e., 
$S(\ell, a) \cup S(\ell', a') = S(\ell'', a'')$ 
for some $\ell''$, $a''$ with $\ell'' > {\rm max} \{\ell,~\ell '\}$. 
It can be easily checked that $S(\ell, a)$, $S(\ell', a')$ are 
adjacent if and only if $a^{-1}a'=q^{\pm i}$ for some 
$$i \in 
\{|\ell -\ell'|+2,\, |\ell -\ell'|+4, \cdots, \ell +\ell' \}.$$ 
Two $q$-strings $S(\ell, a)$, $S(\ell', a')$ are defined to be 
{\it in general position} if they are not adjacent, i.e., 
if either 
\begin{quote}
(i) ~
$S(\ell, a) \cup S(\ell', a')$ is not a $q$-string, \\
or \\
(ii) $\,$
$S(\ell, a) \subseteq S(\ell', a')$ or 
$S(\ell, a) \supseteq S(\ell', a')$. 
\end{quote}
A multi-set ${\{ S(\ell_i,a_i) \}}^n_{i=1}$ 
of $q$-strings 
is said to be {\it in general position} if 
$S(\ell_i,a_i)$ and $S(\ell_j,a_j)$ are in general position for 
any $i,\,j ~(i \neq j,\, 1 \leq i \leq n,\, 1 \leq j \leq n)$. 
The following fact is well-known and easy to prove. 
Let $\Omega$ be a finite multi-set of nonzero scalars from $\C$. 
Then there exists a multi-set 
$\{ S(\ell_i,a_i)\}_{i=1}^n$ 
of $q$-strings in general position such that 
$$ \Omega = \bigcup_{i=1}^n ~S(\ell_i,a_i)$$
as multi-sets of nonzero scalars. 
Moreover such a multi-set of $q$-strings is uniquely determined by $\Omega$. 

\medskip
With a tensor product $V(\ell_1,a_1) \otimes \cdots \otimes V(\ell_n,a_n)$ 
of evaluation modules $V(\ell_i,a_i)~ (1 \leq i \leq n)$, 
we associate the multi-set ${\{ S(\ell_i,a_i) \}}^n_{i=1} $ 
of $q$-strings. 
The following (i), (ii), (iii) are well-known \cite{CP}: 
\begin{quote}
(i)~\,A tensor product 
$V(\ell_1,a_1) \otimes \cdots \otimes V(\ell_n,a_n)$ 
of evaluation modules is irreducible as an $\LL$ -module 
if and only if 
the multi-set ${\{ S(\ell_i,a_i)\}}^n_{i=1}$ of $q$-strings is 
in geneal position. \\
\\
(ii)~\,Set $V=V(\ell_1,a_1) \otimes \cdots \otimes V(\ell_n,a_n)$, 
$V'= V(\ell_1',a_1') \otimes \cdots V(\ell_{n'}',a_{n'}')$ and assume 
that $V$, $V'$ are both irreducible as an $\LL$-module. 
Then $V$, $V'$ are isomorphic as $\LL$-modules if and only if 
the multi-sets ${\{ S(\ell_i,a_i)\}}^n_{i=1}$, 
${\{S(\ell_i',a_i')\}}^{n'}_{i=1}$ 
coincide, i.e., $n=n'$ and $\ell_i=\ell_i'$, $a_i=a_i'$ 
for all $i$ $(1 \leq i \leq n)$ 
with a suitable reordering of 
$S(\ell_1',a_1'), \cdots, S(\ell_n',a_n')$. \\
\\
(iii)~\,Every nontrivial 
finite-dimensional irreducible $\LL$-module of 
type (1,1)  is isomorphic to some 
$V(\ell_1,a_1) \otimes \cdots V(\ell_n,a_n)$.\\
\end{quote}

Two multi-sets 
${\{S(\ell_i,a_i)\}}^n_{i=1}, {\{S(\ell_i',a_i')\}}^{n'}_{i=1}$ 
of $q$-strings are defined to be 
{\it equivalent} if there exists 
$\varepsilon_i \in \{\pm 1\} ~(1 \leq i \leq n)$ 
such that 
${\{S(\ell_i,a_i^{\varepsilon_i})\}}^n_{i=1}$
and 
${\{S(\ell_i',a_i')\}}^{n'}_{i=1}$ 
coincide, i.e., $n=n'$ and 
$\ell_i=\ell_i', \, a_i^{\varepsilon_i}=a_i'$ 
for all $i~(0 \leq i \leq n)$ 
with a suitable reordering of 
$S(\ell_1', a_1'), \cdots, S(\ell_n',a_n')$. 
A multi-set 
${\{S(\ell_i,a_i)\}}^n_{i=1}$
of $q$-strings is defined to be {\it strongly in general position} 
if any multi-set of $q$-strings equivalent to 
${\{S(\ell_i,a_i)\}}^n_{i=1}$ is in general position, i.e., 
the multi-set 
${\{S(\ell_i,a_i^{\varepsilon_i})\}}^n_{i=1}$
is in general position for any choice of 
$\varepsilon_i \in \{ \pm 1 \} ~(1 \leq i \leq n)$.

\begin{lemma}
\label{lemma: q-string}
Let $\Omega$ be a finite multi-set of nonzero scalars from $\C$ 
such that c and $c^{-1}$ appear in $\Omega$ in pairs, i.e., c and 
$c^{-1}$ have the same multiplicity in $\Omega$ for each $c \in \Omega$, 
where we understand that if 1 or -1 appears in $\Omega$, 
it has even multiplicity. 
Then there exists a multi-set 
${\{S(\ell_i,a_i) \}}_{i=1}^n$ of $q$-strings 
strongly in general position such that 
$$ \Omega = \bigcup_{i=1}^n ~
\bigl(S(\ell_i,a_i) \cup S(\ell_i,a_i^{-1})\bigr) $$
as multi-sets of nonzero scalars. 
Such a multi-set of $q$-strings is uniquely determined by $\Omega$ 
up to equivalence. 
\end{lemma}

For an $\LL$-module $V$, let $\rho_V$ denote the representation of $\LL$ 
afforded by the $\LL$-module $V$. Then $\rho_V \circ \varphi_{s,t}$ is 
a representation of $\A$.

\begin{theorem}
\label{thm: q-Onsager}
The following $(i)$, $(ii)$, $(iii)$ hold. 
\begin{enumerate}
\item[$(i)$]
For an $\LL$-module $V=V(\ell_1,a_1) \otimes \cdots \otimes V(\ell_n,a_n)$ 
and nonzero $s,\,t \in \C$, the representaiton 
$\rho_V \circ \varphi_{s,t}$ of $\A$ is irreducible if and only if 
\begin{description}
\item[$(i.1)$] 
the multi-set 
${\{ S(\ell_i,a_i)\}}^n_{i=1}$
of $q$-strings is strongly in general position, 
\item[$(i.2)$] 
none of $-s^2,\, -t^2$ belongs to $S(\ell_i,a_i) \cup S(\ell_i,a_i^{-1})$ 
for any $i$ ~$(1 \leq i \leq n)$, 
\item[$(i.3)$] 
none of the four scalars $\pm st,\, \pm st^{-1}$ equals $q^i$ 
for any $i \in \Z$ ~$(-d+1 \leq i \leq d-1)$. 
\end{description}

\item[$(ii)$]
For $\LL$-modules $V=V(\ell_1,a_1) \otimes \cdots \otimes V(\ell_n,a_n)$, 
$V'=V(\ell_1',a_1') \otimes \cdots \otimes V(\ell_{n'}',a_{n'}')$ 
and $(s,t),\, (s',t') 
\in (\C \backslash \{0\})\times (\C \backslash \{0\})$, 
set $\rho=\rho_V \circ \varphi_{s,t}$ and 
$\rho'=\rho_{V'} \circ \varphi_{s',t'}$. 
Assume that the representations $\rho$, $\rho'$ of $\A$ are both irreducible. 
Then they are isomorphic as representations of $\A$ 
if and only if the multi-sets  
${\{ S(\ell_i,a_i)\}}^n_{i=1}$,  ${\{ S(\ell_i',a_i')\}}^{n'}_{i=1}$ 
are equivalent and 
$$
(s', t') \in 
\{\pm (s,t),\, \pm (t^{-1},s^{-1}), \, \pm (t,s), 
\, \pm (s^{-1},t^{-1})\}.
$$ 

\item[$(iii)$]
Every nontrivial finite-dimensional irreducible representation of $\A$ 
is isomorphic to $\rho_V \circ \varphi_{s,t}$ 
for some $\LL$-module $V=V(\ell_1,a_1) \otimes \cdots \otimes V(\ell_n,a_n)$ 
and $(s,t) \in (\C \backslash \{0\})\times (\C \backslash \{0\})$. 
\end{enumerate}
\end{theorem}

\medskip
Let $A, A^* \in {\rm{End}}(V)$ be a TD-pair \cite{ITT} 
with eigenspaces 
${\{V_i\}}^d_{i=0},~{\{V_i^*\}}^d_{i=0}$ respectively. 
Then we have the split decomposition: 
\[ V= \bigoplus^d_{i=0} U_i , \]
where
\[ U_i=(V_0^* + \cdots + V_i^*) \cap (V_i + \cdots + V_d). \]
By \cite[Corollary 5.7 ]{ITT}, it holds that 
$$ {\rm dim}\, U_i= {\rm dim}\, V_i = {\rm dim}\, V_i^* 
~~~ (0 \leq i \leq d), $$
and 
$$ {\rm dim}\, U_i= {\rm dim}\, U_{d-i} ~~~~~~~~~~
~~~ (0 \leq i \leq d).$$
Note that ${\rm dim}\,U_i$ is invariant under standardization 
of $A,\,A^*$ by affine transformations 
$\lambda A + \mu I$, $\lambda^* A^* + \mu^* I$.  
For an $\LL$-module 
$$V=V(\ell_1,a_1) \otimes \cdots \otimes V(\ell_n,a_n),$$ 
set $A= \rho_V \circ \varphi_{s,t}(z),\,
A^*= \rho_V \circ \varphi_{s,t}(z^*).$ 
If the conditions (i.1), (i.2), (i.3) in Theorem 1 holds, 
then  $A,\,A^*$ are a standardized TD-pair of $q$-Racah type. 
Every standardized TD-pair $A,\,A^*$ of $q$-Racah type arises in this way. 
The split decomposition of $V$ for $A,\,A^*$ coincides 
with the eigenspace decomposition of the element $k_0$ of $\LL$ 
acting on $V$. Thus 
the generating function for 
${\rm dim}\,U_i$
$$ g(\lambda)= \sum^d_{i=0} \,({\rm dim}\, U_i ) \,\lambda^i $$
is given as follows. 

\begin{proposition} 
\label{prop: ch}
{\rm (\cite[Conjecture 13.7 ]{ITT})}
\begin{eqnarray*}
g(\lambda)&=&\prod^n_{i=1}(1+\lambda+\lambda^2+ \cdots + \lambda^{\ell_i}). 
\end{eqnarray*}
\end{proposition}

A TD-pair $A,A^*$ is called a {\it Leonard pair} 
if dim $U_i=1$ for all $i ~~(0 \leq i \leq d)$. 
A standardized TD-pair $A,A^*$ of $q$-Racah type is a Leonard 
pair if and only if it is afforded by an evaluation 
module. In view of this fact, a standardized TD-pair $A,\,A^*$ 
of $q$-Racah type is regarded as a `tensor product of Leonard pairs'.

\medskip
\noindent
Tatsuro Ito \hfil\break
Division of Mathematical and Physical Sciences \hfil\break
Kanazawa University \hfil\break
Kakuma-machi, Kanazawa 920--1192, Japan \hfil\break

\medskip
\noindent
Paul Terwilliger \hfil\break
Department of Mathematics\hfil\break
University of Wisconsin-Madison \hfil\break
Van Vleck Hall \hfil\break
480 Lincoln drive  \hfil\break
Madison, WI 53706-1388, USA \hfil\break


\begin{thebibliography}{8}

\bibitem{CP}
V.~Chari and A.~Pressley, 
\newblock {Q}uantum affine algebras, 
\newblock {\em Commun. Math. Phys.}
{\bf 142} (1991) 261--283. 
\bibitem{ITT}
T.~Ito, K.~Tanabe, and P.~Terwilliger, 
\newblock Some algebra related to ${P}$- and ${Q}$-polynomial association
  schemes,  in:
\newblock {\em Codes and Association Schemes (Piscataway NJ, 1999)}, Amer.
Math. Soc., Providence RI, 2001, 
     167--192; 
{\tt arXiv:math.CO/0406556}.









\bibitem{ITf}
T.~Ito and P.~Terwilliger, 
\newblock The augmented tridiagonal algebra, preprint. 






\end{thebibliography}
\end{document}